\documentclass[12pt]{amsart} 
\usepackage{latexsym}
\usepackage{amssymb}
\usepackage{euscript}

\let\cal\mathcal

\makeatletter \@mparswitchfalse \makeatother

\numberwithin{equation}{section}

\def\ep{{}{\hfill $\Box$} \vskip 5pt \par}

\begin{document}
\setlength{\baselineskip}{18pt}
\title{ Nevanlinna-Pick Kernels and Localization}
\author{Jim Agler}
\thanks{Partially supported by the National Science Foundation}
\author{John E. M\raise.5ex\hbox{\tiny c}Carthy}
\thanks{Partially supported by National Science Foundation
grant DMS 9531967. }
\address{ University of California at San Diego, La Jolla California 92093}
\address{Washington University, St. Louis, Missouri 63130}
\thanks{Math Subject Classification: Primary 47A20, Secondary 46E20}

\bibliographystyle{plain}

\begin{abstract}
We describe those reproducing kernel Hilbert spaces of holomorphic functions
on domains in ${\Bbb C}^d$ for which an analogue of the Nevanlinna-Pick 
theorem holds, in other words when the existence of a (possibly matrix-valued) 
function in the unit ball of the multiplier algebra with specified values 
on a finite set of points is equivalent to the positvity of a related matrix.
Our description is in terms of a certain localization property of the
kernel.
\end{abstract}

\maketitle

\baselineskip = 18pt

\noindent  {\bf \S{0.\ Introduction}}

This paper concerns a generalization of the following result
due to Pick [{\bf P}] and Nevanlinna [{\bf N}].
 \vspace{5mm}
 
\noindent{\bf Theorem 0.1} {\em Let $n$ be a positive integer, let
$\lambda_1, \dots, \lambda_n$ be distinct points in $\Bbb
D$, the open unit disc in the complex plane centered at $0$,
and let $z_1, \dots, z_n \in \Bbb C$.  There exists a
holomorphic function $\varphi$ on $\Bbb D$ with
$\varphi(\lambda_i)=z_i$ for each $i$ and
$\sup\limits_{|z|<1} |\varphi(z)| \le 1$ if and only if the
$n\times n$ matrix $\Big[ (1-z_j\overline z_i)
(1-\lambda_j\overline \lambda_i)^{-1}\Big]$ is positive
semidefinite.}
 \vspace{5mm}

Theorem 0.1 will be generalized in three ways.  First, the
domain $\Bbb D \subseteq \Bbb C$ will be  replaced by a
bounded domain $U$ in $\Bbb C^d$.  Secondly,
observe that if $H^2$ is the  classical Hardy space of
analytic functions on $\Bbb D$ with square integrable
boundary values  and $\varphi$ is a holomorphic function on
$\Bbb D$, then $\varphi$ is a multiplier of $H^2$ (i.e. $\varphi
H^2 \subseteq H^2$) if and only if $\sup\limits_{|z|<1}
|\varphi(z)|$ is finite. Indeed the condition,
$\sup\limits_{|z|<1} |\varphi(z)|\le 1$, in Theorem 0.1 is 
equivalent to the condition that the norm of $\varphi$ as a
multiplier of $H^2 $ (i.e. $\sup\limits_{\stackrel{ f\in H^2} 
{\|f\|=1}} \|\varphi f\|$) is less than or equal to one.
 Furthermore, since the reproducing kernel for
$H^2$ has the form $k_\lambda(\mu) = (1-\overline
\lambda\mu)^{-1}$ it is clear that the $n\times n$ matrix that
appears in Theorem~0.1 can be written in the form
$[(1-z_j\overline z_i) k_{\lambda_i} (\lambda_j)]$.  The second
way in which we shall generalize the statement of Theorem 0.1
will be to replace the bound condition $\sup\limits_{|z|<1}
|\varphi(z)|
\le 1$ on the  holomorphic function $\varphi$ with a bound
condition of the form
$$\sup\limits_{\stackrel{ f\in {\cal H}}{ \|f\|=1}} \|\varphi f\|
\le 1 \leqno(0.2)$$
where $\cal H$ is a Hilbert space of analytic functions on $U$
and to replace the condition that  the $n\times n$ matrix
$[(1-z_j \overline z_i) (1-\overline
\lambda_i\lambda_j)^{-1}]$ be positive semidefinite with
the condition that $[(1-z_j\overline z_i) k_{\lambda_i}
(\lambda_j)]$ be positive semidefinite where
$k_\lambda(\mu)$ is the reproducing kernel for $\cal H$. 
Thirdly, we shall replace the vorgegebene Funktionswerte
$z_1, \dots, z_n$ of Theorem 0.1 with points $z_1, \dots, z_n
\in \cal M_m (\Bbb C)$, the $C^*$-algebra of $m\times m$
matrices.  Since the function $\varphi$ in that event will be
a holomorphic $\cal M_m(\Bbb C)$-valued map on $U$ we
modify (0.2) to the form
$$\sup\limits_{\stackrel{ f\in {\Bbb C}^m\otimes \cal H}{
\|f\|=1}}
\|\varphi f\|
\le 1, \leqno(0.3)
$$ 
and replace the condition that the
$n\times n$ matrix $[(1-z_j\overline z_i) k_{\lambda_i}
(\lambda_j)]$ be positive semidefinite with the condition
that the $mn\times mn$ block matrix $[(1-z_j z^*_i)
k_{\lambda_i}(\lambda_j)]$ be positive  semidefinite.  We
thus have been led to ask the following question.
\vspace{5mm}

\noindent{\bf Question 0.4.} Let $U$ be a bounded domain
 in $\Bbb C^d$ and let $\cal H$ be a Hilbert
space of analytic functions on $\cal H$ with reproducing
kernel $k$.  Fix a positive integer $m$.  When are the
following two conditions equivalent for all choices of
positive integer $n$, all  choices of distinct points
$\lambda_1, \dots, \lambda_n \in U$ and all choices of
points  $z_1,
\dots,z_n \in \cal M_m(\Bbb C)$?
\vspace{3mm}

\begin{enumerate}
\item[(i)] There exists an
$\cal M_m(\Bbb C)$-valued holomorphic function $\varphi$ on
$U$ with
$\varphi(\lambda_i) = z_i$ for each $i$ and with the norm of
$\varphi$ as a multiplier of
$\Bbb C^m \otimes \cal H$ less than or equal to one.

\item[(ii)] The $mn\times mn$ matrix
$[(1-z_jz^*_i) k_{\lambda_i} (\lambda_j)]$ is positive
semidefinite.
\end{enumerate}
\vspace{3mm}

In this paper we shall give an answer to Question 0.4 in the
case when the multiplications by the coordinate functions
of $\Bbb C^d$ form a bounded $d$-tuple of operators $M$
acting  on $\cal H$ with the properties that $\sigma(M) =
cl(U)$ 
and $M^*$ is of sharp type
([{\bf Agr-S}]).  Specifically, let superscripts denote the
standard coordinates in $\Bbb C^d$  and assume for each
integer $r$ in the range $1\le r \le d$ that if $M^r$ is defined
on $\cal H$ by the formula
$$(M^r f) (\lambda) = \lambda^rf(\lambda), f \in \cal H,
\lambda \in U,$$
then
$$M^r\cal H \subseteq \cal H. \leqno(0.5)$$
Evidently, (0.5) will guarantee that, for each $r$, $M^r \in
\cal L(\cal H)$, the $C^*$-algebra of bounded linear
transformations of $\cal H$.  We shall assume that the
$d$-tuple of commuting operators $M=(M^r)$  satisfies the
following conditions.
\vspace{5mm}

\noindent (0.6)\qquad \hangindent=42pt $\sigma(M)$, the 
Taylor spectrum of $M$, (see [{\bf T1}]and [{\bf T2}]) is a
subset of
\newline $cl(U)$ (equivalently, $\sigma(M) = cl(U))$
\vspace{5mm}

\noindent (0.7)\qquad $\sigma_e(M)$, the essential Taylor
spectrum of $M$ (see [{\bf C}]), is a subset of $\partial U$.

$$\dim \bigcap\limits^d_{r=1} \ker (\lambda^r - M^r)^* =1
\quad \text{for\ all}\quad \lambda \in U.
\leqno(0.8)$$
Let us agree to say that a Hilbert space $\cal H$ of
holomorphic functions on a bounded domain
$U$ is {\it regular} if (0.5) - (0.8) are satisfied.
Note that condition (0.8) guarantees that an operator that commutes with
every $M^r$ is multiplication by a  multiplier of $\cal H$, and
conversely all multipliers give rise to operators that commute with $M$
(see {\it e.g.} [{\bf Agr-S}]).

We now describe our answer to Question 0.4.  Let us agree
to say that a Hilbert space  of analytic functions on a
bounded domain is an $m$-interpolation
space if  the two conditions in Question 0.4 are equivalent.  For $\lambda =
(\lambda_1,\dots, \lambda_n)$ an $n$-tuple of distinct
points in $U$, let $\cal H_\lambda = \{f\in \cal H:
f(\lambda_i) =0$ for each $i\}$, and let $\cal S_\lambda$
denote the collection of commuting $d$-tuples of operators
$T$ such that $\sigma(T)\subseteq \{\lambda_1, \dots,
\lambda_n\}$ and $h(T)=0$ whenever  $h$ is holomorphic on
$U$ and $h(\lambda_i)=0$ for each $i$.  For $T$ a commuting
$d$-tuple of operators  let $\cal A_T$ denote the sigma-weak operator topology
closed algebra
generated by the components of $T$.  Finally, let
$H^\infty_k$ denote the algebra of multipliers of $\cal H$. 
Observe that naturally $H^\infty_k \subseteq \cal L(\cal
H)$, so that $\cal M_m(\Bbb C) \otimes H^\infty_k$ carries a
distinguished norm (that inherited from the tensor product
of $\cal M_m(\Bbb C)$ and $\cal L(\cal H)$ as
$C^*$-algebras).  An elementary exercise is to ascertain
that if $\cal M_m(\Bbb C)\otimes H^\infty_k$ is identified
with the space of matrix multipliers of $\Bbb C^m\otimes
\cal H$, then this distinguished norm is the same as the
norm defined in (0.3).

Now, in the ground-breaking papers [{\bf Arv1}] and
[{\bf Arv2}] Arveson introduced the  notion of
$m$-contractivity.  A linear map $\rho$ defined on a
subspace $S$ of a
$C^*$-algebra $A$ and  taking values in a $C^*$-algebra $B$
is said to be $m$-contractive if the map
$$id_m \otimes \rho: \cal M_m(\Bbb C) \otimes
S\longrightarrow \cal M_m(\Bbb C)\otimes B$$
is contractive.  Here, $id_m$ denotes the identity mapping
on $\cal M_m(\Bbb C)$.  A map $\rho$ is said to be
completely contractive if $P$ is $m$-contractive for every
$m$.  Now, let $\cal C_m$ denote the category with objects
the subalgebras of $C^*$-algebras and morphisms the
$m$-contractive algebra homorphisms.  We assume that all
algebras contain a unit and that morphisms  map units to
units.

If $\lambda_1, \dots, \lambda_n$ are $n$ distinct points in
$U$ and $T\in \cal S_\lambda$ observe that $\Phi_T$ defined
via the functional calculus [{\bf T2}] by
$$\Phi_T(\varphi) = \varphi(T),\quad \varphi \in H^\infty_k$$
is an algebra homorphism from $H^\infty_k$ onto $\cal
A_T$.  In particular, if $T_\lambda$ is the commuting
$d$-tuple of operators on $\cal H^\perp_\lambda$ defined by
setting $T^r_\lambda = PM^r|\cal H^\perp_\lambda$ where
$P$ denotes the orthogonal  projection of $\cal H$ onto
$\cal H^\perp_\lambda,\ \Phi_{T_\lambda}$ is a
well-defined algebra homomorphism from $H^\infty_k$
onto $\cal A_{T_\lambda}$.  Also, it is clear that the map
$\Phi_{T_\lambda}$ is completely contractive.  
We shall in future call the  map $\Phi_{T_\lambda}$ the localization
operator and denote it simply $\rho$.  Observe that the
localization operator depends only on $\cal H$ and the
choice of points $\lambda_1, \dots, \lambda_n$.  Let us
agree to say  that $\cal H$ possesses the $m$-contractive
localization property if the following diagram can be
completed in the category $\cal C_m$ whenever $n$ is a
positive integer, $\lambda_1, \dots, \lambda_n$ are $n$
distinct  points in $U$, $T\in \cal S_\lambda$, and $\Phi_T$ is an
$m$-contraction.

\bigskip

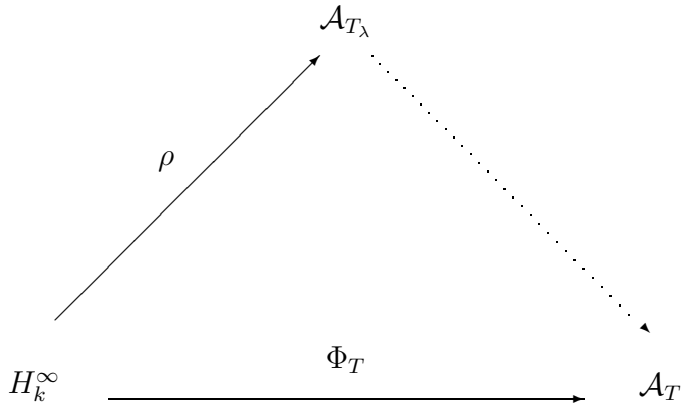
\begin{center}
\begin{picture}(100,200)(100,60)
\put(40,150){\vector(1,1){100}}
\put(60,120){\vector(1,-0){180}}
\put(260,150){\vector(1,-1){05}}
\qbezier[30](160,250)(205,205)(257,152)
\put(148,263){\makebox(0,0){ ${\cal A}_{T_\lambda}$}}
\put(80,210){\makebox(0,0) { $\rho$}}
\put(148,135){\makebox(0,0){ $\Phi_T$}}
\put(30,125){\makebox(0,0) { $H^\infty_k$}}
\put(267,125){\makebox(0,0){ ${\cal A}_T$}}
\end{picture}
\end{center}

\vspace{-2mm}
\noindent Evidently, since $T\in {\cal S}_\lambda, \varphi(T)=0$ 
whenever $\varphi\in \{\varphi \in H^\infty_k:\varphi(\lambda_i)=0\
{\rm for\ each}\ i\} = \ker \rho$, so that the map $\Psi_T: {\cal
A}_{T_\lambda} \to {\cal A}_T$ defined by $\Psi_T(\rho(\varphi))
= \varphi(T)$ always completes the diagram.  Thus, ${\cal H}$ has
the $m$-contractive localization property if and only if
whenever $n$ is a positive integer, $\lambda_1, \dots,
\lambda_n$ are $n$-distinct points in $U, T\in {\cal S}_\lambda$,
and $\Phi_T$ is $m$-contractive,  it is the case that $\Psi_T$ is
also $m$-contractive.
\vspace{2mm}

We now can state our answer to Question 0.4.
\vspace{5mm}

\noindent{\bf Theorem 3.6.}\ {\em Let ${\cal H}$ be a regular
space on a bounded domain $U \subseteq \Bbb
C^d$. 
${\cal H}$ is an $m$-interpolation space if and only if ${\cal H}$
has the
$m$-contractive localization property.}
\vspace{2mm}

We prove this theorem in Section 3.
To understand the basic ideas behind our proof,
we recall some ideas  from [{\bf A1}].   In [{\bf A1}] the collection
of operators that can be modelled using the space ${\cal H}$ was
studied in the case where $d=1$ and $U=\Bbb D$.  Specifically,
one lets $k$ be a positive  definite holomorphic kernel over
$\Bbb D$  with the property that the associated Hilbert space
${\cal H}$ (i.e. the unique Hilbert space of analytic functions with
reproducing kernel $k$) is regular and introduces the collection
of operators ${\cal F}(k)$, consisting of all operators $T$ for
which  there exist a Hilbert space ${\cal K}$, a unital
representation $\pi:{\cal L}({\cal H}) \to {\cal L}({\cal K})$, and a
subspace ${\cal N} \subseteq {\cal K}$ with the properties that
${\cal N}$ is invariant for $\pi(M^*)$ and $T$ is unitarily
equivalent  to $\pi(M^*)|{\cal N}$.  In the case when $k$ is the
Szeg\"o kernel,  $k(\lambda, \mu) =(1-\overline
\lambda \mu)^{-1}, {\cal F}(k)$ is the  set of operators that can be
written as a coisometry restricted to an invariant subspace, the
object of the Sz.-Nagy-Foias and de Branges-Rovnyak model
theories.

There are no surprises in the generalization to several
variables of the setup described in the previous paragraph.  If
$U$ is a bounded domain in $\Bbb C^d$ and ${\cal
H}$ is a  regular Hilbert space of analytic functions on $U$ with
reproducing kernel $k$, then ${\cal F}(k)$ is  defined as the
collection of all commuting $d$-tuples of operators $T=(T^r)$
for which there  exists a Hilbert space ${\cal K}$, a unital
representation $\pi:{\cal L}({\cal H}) \to {\cal L}({\cal K})$, and a
subspace ${\cal N} \subseteq {\cal K}$  with the properties that
${\cal N}$ is invariant for $h(\pi(M^*))$ whenever $h$ is
holomorphic on a neighborhood of $\overline {cl(U)}$ and $T$ is
unitarily equivalent to $\pi(M^*)|{\cal N}$ (i.e. there exists a
unitary operator $X$ such that
$T^r=X^*(\pi(M^r)^*|{\cal N})X$ for each $r$).  
Throughout the paper, we shall use $cl(E)$ to denote the closure of
$E$, and $\overline{E}$ to denote the complex conjugate of $E$.
Note in the 
definition of ${\cal F}(k)$ that 
since $\pi$ is a unital representation,
$\sigma(\pi(M^*)) \subseteq \overline{cl(U)}$ and consequently,
$h(\pi(M^*))$ is well-defined by Taylor's functional calculus
[{\bf T2}].

From now on, we shall set $K=cl(U)$, and let $H(K)$ denote the
space of germs of holomorphic functions on a
neighborhood of $K$.  If $T=(T^r)$ is a commuting $d$-tuple of
operators acting on a space ${\cal K}$ and 
$\sigma(T)\subseteq K$, we say a subspace ${\cal N} \subseteq
{\cal K}$ is $H(K)$-invariant for $T$ if $h(T){\cal N} \subseteq
{\cal N}$ whenever $h \in H(K)$.  Thus, if $T$ is a commuting
$d$-tuple with $\sigma(T) \subseteq {\overline K}, T \in {\cal F}(k)$ if and
only if there exist a unital representation $\pi$ of ${\cal L}({\cal
H})$ and an $H({\overline K})$-invariant subspace ${\cal N}$ for $\pi(M^*)$
such that $T$ is unitarily equivalent to $\pi(M^*)|_{\cal N}$.  We shall
abbreviate this latter condition in the sequel and simply say
that $T$  has an $H({\overline K})$-extension to a tuple of the form
$\pi(M^*)$.

Now a fundamental fact about the family ${\cal F}(k)$ is that
$H({\overline K})$-extensions localize.  Precisely what this means is that if
$\lambda_1, \dots, \lambda_n$ are $n$ distinct points in $U, T \in
{\cal S}_{\overline \lambda}$, and $T$  has an $H({\overline K})$-extension to
a tuple of the form $\pi(M^*)$, then $T$ has an extension to a 
tuple of the form $\pi(M^*_z | {\cal H}^\perp_\lambda)$.  This
result, which is a key element in the proof of Theorem 3.6, is
proved in Section 1 of this paper (Theorem 1.2).

Since extensions localize it is natural to ask whether dilations
localize.  Let us agree  to say that a commuting $d$-tuple
$T=(T^r)$ has an $H({\overline K})$-dilation to a tuple of the form
$\pi(M^*_z)$ if there exist a unital representation $\pi$ of ${\cal
L}({\cal H})$ and a pair of $H({\overline K})$-invariant subspaces ${\cal N}_1,
{\cal N}_2$ for $\pi(M^*)$ with the properties that ${\cal N}_1
\subseteq {\cal N}_2$ and $T$ is unitarily equivalent  to 
$P_{{\cal
N}_2 \ominus {\cal N}_1}\pi(M^*)|{\cal N}_2 \ominus {\cal N}_1$. 
Let us agree to say ${\cal H}$ has the {\it dilation localization
property}  if whenever $n$ is a positive integer, $\lambda_1,
\dots, \lambda_n$ are $n$ distinct points in $U, T\in {\cal
S}_{ \lambda}$, and $T^*$ has an $H({\overline K})$-dilation to an
operator of the form $\pi(M^*)$, then $T^*$ has an $H({\overline K})$-dilation
to  an operator of the form $\pi(M^* | {\cal H}^\perp_\lambda)$.
Equivalently, if $T$ has an $H(K)$-dilation to an operator of the form 
$\pi(M)$, then $T$ has an $H(K)$-dilation to an operator of the form
$\pi(PMP)$, where $P$ is the projection from $\cal H$ onto 
${\cal H}^\perp_\lambda)$.

Now a straightforward consequence of 
Arveson's
theory of completely contractive algebra homomorphisms is
that ${\cal H}$ has the dilation localization property if and only if
${\cal H}$ has the $m$-contractive localization property for
every $m$.  We thus obtain from Theorem 3.6 the following
result:

\vspace{5mm}

\noindent{\bf Theorem 3.5.}\ \ {\em Let ${\cal H}$ be a regular
space on a bounded domain $U\subseteq \Bbb
C^d$. 
${\cal H}$ is a complete interpolation space if and only if ${\cal
H}$ has the $H({\overline K})$-dilation localization property.}
\vspace{2mm}

In Section 4 we give two concrete applications of Theorem 3.5;
one in the multiplier norm of the Dirichlet space in one variable
and the other in a norm on a space of holomorphic functions on the ball in
$\Bbb C^d$ that agrees with the ordinary $H^\infty$ norm when
$d=1$.  These  applications demonstrate that Theorem 3.6 and
3.5 already contain whatever concrete function theory one
might believe is involved in the equivalence of 0.4 (i) and 0.4
(ii).

In this paper we shall always adhere to the following notations. 
$T$ will denote  a commuting $d$-tuple of bounded
operators acting on a Hilbert space.  The components  of $T$ (as
well as the components of elements of $\Bbb C^d$) will be
denoted with superscripts.  All operators on $d$-tuples will be
assumed to act componentwise.  Thus, for example, if $T$ is a
commuting $d$-tuple of operators acting on a space ${\cal H}$,
then $T^* = (T^{r^*})$ and if $\pi:{\cal L}({\cal H}) \to {\cal L}({\cal
K})$ is a representation then $\pi(T) = (\pi(T^r))$.  If
$\{T_\alpha:\alpha \in {\frak A}\}$ is a collection of $d$-tuples
with $T_\alpha$ acting on ${\cal H}_\alpha$, then
$\bigoplus\limits_{\alpha\in{\cal H}_\alpha} T_\alpha$ denotes
the $d$-tuple $\left(\bigoplus\limits_{\alpha\in{\frak A}}
T^r_\alpha\right)$ which acts on
$\bigoplus\limits_{\alpha\in{\frak A}} {\cal H}_\alpha$.  If ${\cal
H}$ is a Hilbert space and ${\cal S} \subseteq {\cal H}$, then
$[s:s\in S]$ denotes the closed linear span of $S$ in ${\cal H}$ and
if $x, y\in {\cal H}$, then $x\otimes y \in {\cal L}({\cal H})$ is
defined by $x\otimes y(h) = \langle h, y\rangle x$.

Finally, we remark that because of the fundamental work of
Aronszajn [{\bf Aro}], it is well  known that ${\cal H}$ and $k$
determine each other uniquely.  Accordingly, if a kernel is
assumed to be regular, this means that the corresponding
Hilbert space $\cal H$ is regular.  Similarly, we  make no
particular distinction between kernels and the spaces
determined by them with respect to the notions of being
$m$-contractively localizable, dilation localizable, having the
$m$-interpolation property, or having the complete
interpolation property.

\vspace{1cm}

\noindent {\bf  \S{1.\ 
Coanalytic Models}}
 \vspace{2mm}

In this section $U$ will be a bounded domain in
$\Bbb C^d$, we shall set $K=cl(U)$ and ${\cal H}$ will be a regular space
on $U$ with reproducing kernel $k$.  
We emphasize that we are assuming $\sigma(M) \subseteq K$ and
$\sigma_e(M) \subseteq \partial U$.  Let $\cal F$ denote the
model generated by $M$, i.e. the collection of all commuting
$d$-tuples of operators on Hilbert space that are of the form
$\pi(M^*)|{\cal N}$ where $\pi$ is a unital  representation of ${\cal
L}({\cal H})$ and $\cal N$ is an $H(\overline{K})$-invariant subspace for
$\pi(M^*)$.  Our first result is an obvious extension of Theorem
2.8 in [{\bf A1}].  Accordingly, we merely provide a brief sketch
of its proof.
\vspace{5mm}

\noindent{\bf Theorem 1.1.}\ \ {\em Let $\cal K$ be a Hilbert
space and let $J\in {\cal L}({\cal K})^{(d)}$.  Let ${\cal A}$ be a 
$C^*$-algebra containing $M$ and all compact operators on ${\cal H}$.
There exists a unital
representation $\pi:{\cal A} \to {\cal L}({\cal K})$ with $\pi(M)=J$ if
and only if $J$ is unitarily equivalent to a tuple of the form
$M^{(\nu)}\oplus \pi_0(M)$ where $\nu$ is a cardinal and $\pi_0$
is a unital representation of ${\cal A}$ that annihilates
the compact operators on ${\cal H}$.}
\vspace{5mm}

\noindent{\bf Proof.}\  First assume that $J\in {\cal L}({\cal K})^{(d)},
{\cal K}_0$ is a Hilbert space, $\pi_0:{\cal A} \to {\cal
L}({\cal K}_0)$ is a  unital representation killing the compacts,
and
$U:{\cal K} \to H^{(\nu)} \oplus {\cal K}_0$ is a Hilbert space
isomorphism with $J=U^*(M^{(\nu)} \oplus \pi_0(M))U$.  Define
$\pi:{\cal A} \to {\cal L}({\cal K})$ by 
$$\pi(X) = U^*(X^{(\nu)}\oplus \pi_0(X))U.$$
Obviously, $\pi$ is a unital representation and $\pi(M)=J$.
\vspace{2mm}

Conversely, assume that $\pi:{\cal A} \to {\cal L}({\cal
K})$ is a unital representation and \newline $J=\pi(M)$.  Setting
$${\cal K}_1 = [\pi(f\otimes g)x: f, g \in {\cal H},\  x\in {\cal K}]$$
and
$${\cal K}_0 = {\cal K} \ominus {\cal K}_1,$$
it is easy to verify that $\pi=\pi_1\oplus \pi_0$ where
$\pi_1:{\cal A} \to {\cal L}({\cal K}_1)$ and $\pi_0:
{\cal A} \to {\cal L}({\cal K}_0)$ are unital representations.  An
analysis of the definition of ${\cal K}_1$ reveals that $\pi_0$
kills the  compact operators on ${\cal H}$.  Finally, observing that 
$${\cal K}_1 = [\pi(k_\lambda\otimes g)x: \lambda \in U, g\in
{\cal H}, x\in {\cal K}]$$
allows one to deduce as in [{\bf A1}] that $\pi_1(M)$ is unitarily
equivalent to $M^{(\nu)}$ where 
$$\nu = \dim([\pi(k_\lambda \otimes g)x: g \in {\cal H}, x\in {\cal
K}])$$
is independent of the choice of $\lambda$. This establishes
Theorem 1.1.
\ep

In the project of studying the contractions model theoretically
the von Neumann-Wold decomposition theorem plays a crucial
role.  Once one has realized that it is the  coisometries that will
form the collection of modelling operators (i.e. the operators to
which the general contraction extends) then the fact that the
general coisometry has a  particular concrete form is an
important step in this project. In the present context where  the
contractions have been replaced by the more general class
$\cal F$ the von Neumann-Wold decomposition is replaced by
Theorem 1.1.

Fix $n$
distinct points $\lambda_1, \dots, \lambda_n\in U$ and consider
the ideal $I_{\overline \lambda}\subseteq \Bbb C[x_1, \dots, x_d]$, the ring of
polynomials in $d$ variables, defined by 
$$I_{\overline \lambda} = \{p: p(\overline \lambda_i) =0  \quad
{\rm whenever}\quad 1\le i \le n\}.$$
Associated with $I_{\overline \lambda}$ is the localized model
${\cal F}_{\overline \lambda}$ which is defined as the set of
all $T\in {\cal F}$ such that  $p(T)=0$ whenever $p\in I_{\overline
\lambda}$.  The following localization result is the key to
studying  the norm in the $n$-dimensional Banach algebra
formed from $H^\infty_k$ by factoring out
$(H^\infty_k)_\lambda$, the ideal of functions in $H^\infty_k$
that vanish at the points $\lambda_1, \dots, \lambda_n$.

\vspace{5mm}

\noindent{\bf Theorem 1.2.}\ \ {\em Let $\cal G$ be a Hilbert
space and let $T\in {\cal L}({\cal G})$.  $T\in {\cal F}_{\overline 
\lambda}$ if and only if there exists a cardinal $\nu$ 
and an invariant subspace
$\cal N$ for $(M^* |\,[k_{\lambda_1},
\dots,k_{\lambda_n}])^{(\nu)}$ such that $T$ is unitarily
equivalent to $(M^* |\, [k_{\lambda_1},
\dots k_{\lambda_n}])^{(\nu)} | {\cal N}$.}

\vspace{5mm}

\noindent{\bf Proof.}\ \ First observe that $M^*| [k_{\lambda_1},
\dots k_{\lambda_n}] \in {\cal F}_{\overline \lambda}$ and
that ${\cal F}_{\overline \lambda}$ is closed with respect  to
the operations of forming direct sums and restricting to
invariant subspaces.

Conversely assume that $T\in {\cal F}_{\overline \lambda}$.
Theorem 1.1 and the definition of $\cal F$ imply that there
exists a cardinal $\nu$ and a unital representation $\pi_0 : {\cal
L}(H^2_i)\to \cal L(\cal K_0)$ which annihilates the compact
operators on ${\cal H}$ and an isometry $V:\cal G\to
({\cal H})^{(\nu)} \oplus \cal K_0$ such that 
$$
\begin{array}{l}
h(T) = V^*[h(M^*)^{(\nu)} \oplus h(\pi_0(M^*))]V\\
{\rm and} \\
VV^* [ h(M^*)^{(\nu)} \oplus h(\pi_0(M^*)) ] VV^* =
[ h(M^*)^{(\nu)} \oplus h(\pi_0(M^*)) ] VV^*
\end{array}
\leqno(1.3)$$
whenever $h\in H(\overline K)$.

We claim first that in (1.3) it may be assumed that the
$h(\pi_0(M^*))$ summand is absent.  Equivalently, if $V$ is
decomposed,
$$V= \left[
\begin{array}{c}
V_1\\ V_0
\end{array}\right],$$
with respect to $({\cal H})^{(\nu)}\oplus \{0\}$ and $\{0\}\oplus
\cal K_0$, then $V_0 = 0$.  To see this recall that
$\sigma_e(M)\subseteq \partial U$  and $\pi_0$ kills the
compacts.  Hence $\sigma(\pi_0(M))\subseteq \partial U$.  On
the other hand (1.3) implies that 
$$V^*_0p(\pi_0(M^*))^* p(\pi_0(M^*))V_0 =0$$
whenever $p\in I_{\overline \lambda}$. But these facts would
imply a contradiction if $V_0\not= 0$.  For suppose $x\in {\cal
G}$ and $V_0 x \not=0$.  Let $i$ denote the first positive integer
such that
$$p(\pi_0(M^*))V_0 x = 0$$
whenever $p\in I_{\{ \overline \lambda_1, \dots, \overline
\lambda_i\}}$ and choose $p_0\in I_{\{\overline \lambda_1, \dots,
\overline
\lambda_{i-1}\}}$ such that
$$p_0(\pi_0(M^*))V_0 x \not= 0.$$
By construction,
$$p_0(\pi_0(M^*))V_0 x\in \bigcap\limits^d_{r=1} \ker
\left(\overline{\lambda^{r}_i} - \pi_0 (M^r)^*\right),$$
that is
$$\overline \lambda_i \in \sigma_p(\pi_0(M^*)).$$  
This contradicts $\sigma(\pi_0(M))\subseteq \partial U$  and
establishes our claim that the $\pi_0(M^*)$ summand is absent.

We have shown that (1.3) can be reformulated as:
$$h(T) = V^*h(M^*)^{(\nu)}V\leqno(1.4)$$
whenever $h\in H(\overline K)$.  Since $p(T)^*p(T)=0$ whenever
$p\in I_{\overline \lambda}$, it follows from (0.8) that
ran$V\subseteq [k_{\lambda_1}, \dots, k_{\lambda_n}]^{(\nu)}$
which concludes the proof of Theorem 1.2.
\ep

A basic
fact in Taylor's functional calculus [{\bf T2}] is that if $T$ is a 
commuting $d$-tuple of operators and $\sigma(T)\subseteq K$
then the map
$$h\longmapsto h(T)$$
defined on $H(K)$ is continuous on $H(K)$.  If $f\in H(K)$ define
$\check f\in H(\overline K)$ by setting $\check f(\lambda) =
\overline{f(\overline \lambda)}$.  It is then the case that
$\check f(T^*)  =f(T)^*$ whenever $f\in H(K)$ and
$\sigma(T)\subseteq K$.  Also the map 
$$\check h \longmapsto \check h(T^*)$$
 is continuous on $H(\overline K)$.  It follows by the nuclearity 
of $H(K)$ that the bilinear map defined  on $H(\overline K)\times
H(K)$ by 
$$(\check g, f) \longmapsto \check g(T^*) f(T)$$
can be extended uniquely to a continuous map
$$h\longmapsto h(T)$$
defined on $H(\overline K\times K)$  hereafter referred to as
the hereditary functional calculus for $T$.  The following result
is then clear (Theorem 1.5 in [{\bf A1}]).

\vspace{5mm}

\noindent{\bf Theorem  1.3.}\ \ {\em Let $U$ be a bounded domain and let $\cal H$ be a regular space  on $U$ with
associated kernel $k$.  $T\in {\cal F}(k)$ if and only if $\sigma(T)
\subseteq \overline K$ and the map 
$$h(M^*) \longmapsto h(T)\ ,\quad  h\in H(K \times \overline K)$$
is completely positive.}

\vspace{2mm}

As in [{\bf A1}] Theorem 1.3 will allow us to derive a concrete
condition for a tuple $T$ with  $\sigma(T) \subseteq \overline U$
to be an element of ${\cal F}(k)$.  The key is the calculation that
proves the following lemma whose proof is identical to the
proof of Proposition 2.5 in [{\bf A1}].

\vspace{5mm}

\noindent{\bf Lemma 1.4.}\ \ {\em Let $U$ be a bounded domain and let $\cal H$ be a regular space over $U$ with
kernel $k$.  If $h\in {\cal M}_m(\Bbb C) \otimes H(K\times
\overline K)$, then $h(M^*) \ge 0$ if and only if the ${\cal M}_m
(\Bbb C)$-valued holomorphic kernel function,
$$h(\mu, \overline \lambda) k_\lambda(\mu)$$ 
is positive semi-definite on $U$.}

We conclude this section with the promised concrete condition for
$T\in {\cal F}(k)$ when $\sigma(T)\subseteq U$.  For its proof
follow the proof of Theorem 2.3 in [{\bf A1}]. Remember that $\overline U$
is the complex conjugate of $U$.

\vspace{5mm}

\noindent{\bf Theorem 1.5.}\ \ {\em Let $U$ be a bounded domain in $\Bbb C^d$ and let $\cal H$ be a regular space
over $U$ with nonvanishing kernel $k$.  If $T$ is a commuting
$d$-tuple of operators with $\sigma(T) \subseteq \overline U$,
then $T\in {\cal F}(k)$ if and only if $\displaystyle\frac{1}{k} (T)
\ge 0.$}

\vspace{1cm}

\noindent {\bf \S2. Discrete Matrix Interpolation}

In this section we shall review some of the ideas in Section I of
[{\bf A2}] and then derive a generalization  of Proposition 1.18
of that paper.  Our exposition will be somewhat terse; the
reader is invited to consult [{\bf A2}] for a chatty discussion.

Fix a positive integer  $n$ and let $N=\{i \in \Bbb Z: 1 \le i \le
n\}$.  If $g:N\times N\to \Bbb C$ we say
$g\ge 0$ if the $n\times n$ matrix $(g(i, j))_{i, j\in N}$ is positive
semidefinite.  We define the support of $g, spt(g)$, by
$$spt(g) = \{i \in N: g(i, i) \neq 0\}.$$
If $g:N\times N\to \Bbb C$ and $g\ge 0$ we say $g$ is a kernel on
$N$ if the matrix
$$(g(i, j))_{i, j\in spt(g)}$$
is positive definite.  If $g$ is a kernel on $N$, we define for each
$i \in N, g_i: N\to \Bbb C$, by the  formula $g_i(j) = g(i, j)$.  We
then form a Hilbert space $H_g$ by defining an inner product on 
linear combinations of the form $\sum\limits_{i \in spt(g)} c_i
g_i$ by setting
$$\langle \sum\limits_i c_ig_i,\  \sum\limits_j d_j g_j \rangle =
 \sum\limits_{i, j} g(i, j) c_i \overline d_j.$$
Observe that $H_g$ is a Hilbert space of functions on $N$ and
that $g_i$ has the reproducing  property $\langle f, g_i\rangle
=f(i)$ whenever $f\in H_g$.

The alert reader will have noticed that we have just duplicated
the construction that  this paper began with, with the set $U$
replaced by the set $N$ but without the hypothesis  that $g$ be
strictly positive definite.  The analog of the operators $h(M)^*$
for $h\in H(K)$ would now be the class of operators with the
property that the $g_i$ are eigenfunctions for  the operators. 
Accordingly for $z:N\to \Bbb C$ define $T_{g, z} \in {\cal L}(H_g)$
by requiring that 
$$T_{g, z} g_i = \overline {z(i)} g_i\ ,\ \ i\in\ spt(g).$$
The adjoint is taken in this formula so that the notation will be
consistent with (0.5).  Also observe that $T_{g, z}$ depends only
on the values of $z$ on the set $spt(g)$.

Now if $I\subseteq N$ define $H_g(I) \subseteq H_g$ by $H_g(I)
= [g_i: i\in I\cap spt(g)]$.  If $I\subseteq N$ and  $z:I\to \Bbb C$,
define $T_{g, z}\in {\cal L}(H_g(I))$ by requiring that 
$$T_{g, z} g_i = \overline {z(i)} g_i\ \ i \in I \cap spt(g).$$
If one knows $z$ then one knows the domain of $z$ so this
notation is unambiguous.  Also observe that if $I_0 \subseteq
I_1 \subseteq N,\  z_1:I_1 \to \Bbb C$ and $z_0=z_1 |I_0$, then
$H_g(I_0)$ is invariant for  $T_{g, z_1}$ and 
$$T_{g, z_0} = T_{g, z_1} | H_g(I_0).$$
Now, if $I_0 \subseteq I_1 \subseteq N$ and $z:I_0 \to \Bbb C$,
define $H_g(I_0, I_1) \subseteq H_g$ by 
$$H_g(I_0, I_1) = H_g(I_1) \ominus H_g(I_1\backslash I_0)$$
and define $T_{g,z}(I_1) \in {\cal L}(H_g(I_0, I_1))$ by 
$$T_{g, z}(I_1) = PT_{g,\tilde z} | H_g(I_0, I_1),$$
where $P$ is the orthogonal projection of $H_g(I_1)$ onto
$H_g(I_0, I_1)$ and $\tilde z$ is any extension of  $z$ to $I_1$
(i.e. $\tilde z: I_1 \to \Bbb C$ and $\tilde z |I_0 =z)$.  It can be
shown that $T_{g, z}(I_1)$ depends only on $z$  and $I_1$.  It does
not depend on the choice of extension $\tilde z$.

 We now extend these ideas in an obvious way to the vector
valued case.  For $m$ a  positive integer let $H_{g, m} = \Bbb
C^m \otimes H_g$.  Denoting $a \otimes g_i \in H_{g, m}$ by
simply $ag_i$ when $a \in \Bbb C^m$  and $i \in spt(g)$ it is clear
that the general element $f\in H_{g, m}$ can be represented
uniquely
in the form $\sum\limits_{i \in spt(g)} a_ig_i$.  If $I\subseteq N$
set $H_{g, m}(I) = [ag_i: a\in \Bbb C^m, i \in spt(g) \cap I]$.  If
$I\subseteq N$ and $z:I\to {\cal M}_m(\Bbb C)$ define $T_{g, z}
\in {\cal L}(H_{g, m}(I))$ by requiring that 
$$T_{g, z}(ag_i) = (z(i)^*a)g_i,\ a\in \Bbb C^m,\ i\in spt(g) \cap
I.$$
As before observe that if $I_0\subseteq I_1 \subseteq N, z_1:
I_1 \to {\cal M}_m(\Bbb C)$, and $z_0=z_1|I_0$, then 
$H_{g, m}(I_0)$ is invariant for 
$T_{g, z_1}$ and 
$T_{g,z_0} =T_{g, z_1}
|H_{g, m}(I_0)$. If $I_0\subseteq I_1 \subseteq N$, define
$H_{g,m}(I_0, I_1) = H_{g, m}(I_1) \ominus H_{g,m}(I_1\backslash
I_0)$.  If $I_0\subseteq I_1 \subseteq N$ and $z: I_0\to {\cal
M}_m(\Bbb C)$ define $T_{g, z}(I_1) \in {\cal L}(H_{g, m}(I_0, I_1))$ by
$$T_{g, z}(I_1) = PT_{g, \tilde z}| H_{g, m} (I_0, I_1),$$
where $P$ is the orthogonal projection of $H_{g, m}(I_1)$ onto
$H_{g, m}(I_0, I_1)$ and $\tilde z$ is any extension of $z$ to
$I_1$.  As before the definition of $T_{g, z}(I_1)$ does not
depend on the choice of $\tilde  z$.

We now are ready to state and prove the promised
generalization of Proposition 1.18 from [{\bf A2}].  Let us agree
to say $g$ is a discrete $m$-interpolation kernel on $N$ if for all
$I\subseteq N$  and all $z:I\to \Bbb C$
$$\|T_{g, z}\| = \inf\limits_{\stackrel{\tilde z:N\to {\cal
M}_m(\Bbb C)}{\tilde z|_I=z}}\ \|T_{g, \tilde z}\|.$$
\vspace{5mm}

\noindent{\bf Proposition 2.1.}\ \ {\em Let $g$ be a kernel on
$N$.  The following three conditions are equivalent.
\begin{enumerate}
\item [(i)]$g$ is a discrete $m$-interpolation kernel.

\item[(ii)]  If $I_0\subseteq I_1 \subseteq N$ and $z:I_0\to {\cal
M}_m(\Bbb C)$, then
$$\|T_{g, z}(I_1) \| \le \|T_{g, z}\|.$$
\item[(iii)] If $I_0\subseteq I_1 \subseteq N,\ I_1\backslash
I_0$  consists of a single point and $z:I_0 \to {\cal M}_m(\Bbb
C)$, then 
$$\|T_{g, z}(I_1) \| \le \|T_{g, z}\|.$$
\end{enumerate}
}

\vspace{5mm}

\noindent{\bf Proof.}\ \  That (i) $\Rightarrow$ (ii) follows in
exactly the same way as in the proof of Proposition 1.18 in [{\bf
A2}].  
Indeed, by hypothesis (i), any $z$ defined on $I_0$
can be extended to $\tilde z$ defined on all of $N$
so that $\|T_{g, \tilde z}\|$ is arbitrarily close to 
$\|T_{g, z}\|$; and compressing $T_{g, \tilde z}$ to 
$H_g(I_0, I_1)$ can not increase the norm.

Obviously (ii) $\Rightarrow$ (iii).  There remains to show
that (iii) $\Rightarrow$ (i).   Accordingly fix a kernel $g$,
assume that (iii) holds, let $I\subseteq N, z:I\to {\cal M}_m(\Bbb
C)$, let
$i' \in N\backslash I$ and set  $I'=I\cup\{i'\}$.  We shall show that
$$\inf\limits_{\stackrel{w:I'\to {\cal M}_m(\Bbb C)}{w|_I=z}}\ 
\|T_{g, w}\| =\|T_{g, z}\|.\leqno(2.2)$$
Since $I, i'$, and $z$ are arbitrary, condition (i) will then follow by
iteration.

Now by an argument similar to that which occurs in[{\bf A2}], the
infimum in (2.2) is actually attained.  Choose $\tilde z:I'\to \Bbb
C$ such that $\tilde z|I=z$ and such that
$$\rho_1=\|T_{g, \tilde z}\|^2 = \inf\limits_{\stackrel{w:I'\to
{\cal M}_m(\Bbb C)}{w|I=z}} \ \|T_{g, w}\|.$$
Set $\rho_0 = \|T_{g,  z}\|^2$.  Thus, Proposition 2.1 will be
established if we can show that $\rho_1=\rho_0$.  We shall
argue by contradiction.  Accordingly assume that 
$$\rho_0 <\rho_1. \leqno(2.3)$$

Now choose $\omega \in H_g(I') \ominus H_g(I)$ with the
property that $\langle g_{i'}, \omega\rangle=1$.  The existence
of $\omega$ is guaranteed by (2.3).  Fix an orthonormal basis
$\{e_r\}\subseteq \Bbb C^m$ and observe that if $\delta \in {\cal
M}_m(\Bbb C)$ then 
$$\varphi_\delta(t) = \|T_{g, \tilde z} + t \sum\limits_r(\delta
e_r) g_{i'} \otimes e_r\omega \|^2$$
 is a candidate for the infimum in (2.2).  Now $\varphi_\delta$ is
not in general differentiable at $t=0$.  However it is
differentiable from the right at $t=0$ and its derivative can be
calculated in the following manner.  Set ${\cal M}=\ker(\rho_1 -
T_{g, \tilde z}^* T_{g, \tilde z})$.  Then
$$\lim\limits_{t\to 0 +}\  \frac{1}{t}\ (\varphi_\delta(t) -
\varphi_\delta(0)) =\sup\limits_{\stackrel{\gamma \in {\cal
M}}{\|\gamma\|=1}}\ 2 Re \langle T_{g, \tilde
z}^*\left(\sum\limits_r(\delta e_r) g_{i'} \otimes
e_r\omega\right) \gamma, \gamma\rangle.\leqno(2.4)$$
Since $0$ is a local minimum for $\varphi_\delta$ whenever
$\delta \in {\cal M}_m(\Bbb C)$ we deduce from (2.4) that 
$$0\le \sup\limits_{\stackrel{\gamma \in {\cal
M}}{\|\gamma\|=1}}\ 2 Re \langle T_{g, \tilde
z}^*\left(\sum\limits_r(\delta e_r) g_{i'} \otimes
e_r\omega\right) \gamma, \gamma\rangle\leqno(2.5)$$
whenever $\delta \in {\cal M}_m(\Bbb C)$.  Now let $\cal P$
denote the convex set of 
positive
operators $A\in {\cal L}(H_{g, m}(I'))$
with the properties ran$A\subseteq {\cal M}$ and $trA=1$.  Since
the operators of the form $\gamma \otimes \gamma$ with
$\gamma \in {\cal M}$ and $\|\gamma\|=1$ are in $\cal P$ we
deduce from (2.5) that
$$0\le \inf\limits_{\delta \in {\cal M}_m(\Bbb C)}\
\sup\limits_{A\in {\cal P}}\ 2
Re\ tr\left(T^*_{g, \tilde z}\left(\sum\limits_r(\delta e_r)
g_{i'}
\otimes e_r\omega\right)A\right).\leqno(2.6)$$
Now observe in (2.6) that the objective in the min-max
problem is a real bilinear function in $\delta$ and $A$. 
Consequently, the von Neumann minimax theorem (or indeed
the Hahn-Banach theorem) implies that there exists $A_0 \in
{\cal P}$ such that 
$$0\le  2 Re\ tr \left( T_{g, \tilde
z}^*\left(\sum\limits_r(\delta e_r) g_{i'} \otimes
e_r\omega\right) A_0 \right)$$
for all $\delta \in {\cal M}_m(\Bbb C)$.  Replacing $\delta$ by
$e^{i\theta}\delta$ with $e^{i\theta}$ appropriately chosen thus
yields that in fact,
$$0= tr \left( T_{g, \tilde
z}^*\left(\sum\limits_r(\delta e_r) g_{i'} \otimes
e_r\omega\right) A_0\right)$$
for all $\delta \in {\cal M}_m(\Bbb C)$.  In this last equality
letting $\delta=x\otimes e_s$ where $x\in \Bbb C^m$ and $1\le
s \le m$ reveals that
$$0=tr\  T_{g, \tilde
z}^*(x g_{i'} \otimes
e_s\omega ) A_0$$
whenever $x\in \Bbb C^m$ and $1 \le s \le m$.  Hence
$$T_{g, \tilde
z}^*(\Bbb C^m g_{i'}) \perp A_0(\Bbb C^m \omega).\leqno(2.7)$$

Now, we claim that $A_0(\Bbb C^m\omega) \not= \{0\}$.  For
otherwise, $\Bbb C^m\omega \perp$ ran$A_0 \subseteq {\cal
M}$.  Since $tr\ A_0 =1$ this would imply the existence of a
nonzero vector in ${\cal M}$ (the space on which  $T_{g, \tilde
z}$ attains its norm) that is orthogonal to $\Bbb C^m\omega$
contradicting (2.3).  Choose $y \in A_0(\Bbb C^m \omega)$ with
$\|y\|=1$.  Observe that since ran$A_0 \subseteq {\cal M}, T_{g, \tilde
z}^* T_{g, \tilde
z} y = \rho_1 y$.  Also (2.7) implies that 
$$T_{g, \tilde
z} y \perp \Bbb C^m g_{i'}.\leqno(2.8)$$
Now, if $a\in \Bbb C^m$, then
\begin{eqnarray*}
\langle y, ag_{i'}\rangle & = & \frac{1}{\rho_1}\ \langle T_{g,
\tilde z}^* T_{g, \tilde
z} y, ag_{i'}\rangle\\
&=& \frac{1}{\rho_1} \langle T_{g, \tilde
z} y, (\tilde z(i')^*a) g_{i'}\rangle \\
&=& 0.
\end{eqnarray*}
Hence we also have that 
$$y \perp \Bbb C^m g_{i'}. \leqno(2.9)$$

We now use (iii), (2.8), and (2.9) to derive a contradiction to
(2.3).  Let $P$ denote the orthogonal projection of ${\cal H}_{g,
m}(I')$ onto ${\cal H}_{g, m}(I, I')$.  Observe that (2.8) and (2.9)
imply that $y$ and $ T_{g, \tilde z} y $ are in 
${\cal H}_{g, m}(I, I')$.  Hence
\begin{eqnarray*}
\rho_1 &=& \|T_{g, \tilde z}y\|^2 \\
& = & \|(PT_{g, \tilde z}| {\cal H}_{g,m} (I, I'))y\|^2\\
&=& \|T_{g, \tilde z}(I')y\|^2 \\
&\le& \|T_{g, \tilde z}(I')y\|^2\\
&\le& \|T_{g, z}\|\\
&=&\rho_0,
\end{eqnarray*}
a contradiction which concludes the proof of Proposition 2.1.
\ep

\vspace{1cm}

\noindent {\bf  \S{3.\ 
Holomorphic Interpolation Kernels}}
 \vspace{2mm}

In this section we shall give a concrete model theoretic
condition on the family ${\cal F}(k)$ which is both necessary and
sufficient for the Nevanlinna-Pick interpolation result to be
true for the space $H^\infty_k$.  Specifically, we introduce the
following definition.

\vspace{5mm}
\noindent{\bf Definition 3.1.}\ \   Let $U\subseteq \Bbb C^d$ be a
bounded domain and let $k$ be a kernel on $U$ of
the type described in the introduction. Then $k$ is a holomorphic
$m$-interpolation kernel on $U$ if for each positive integer $n$,
each choice of distinct points $\lambda_1, \dots, \lambda_n\in
U$ and each choice of matrices $z_1, \dots, z_n \in {\cal
M}_m(\Bbb C)$, there exists a function $\varphi \in H^\infty_{k,
m}$ with $\|\varphi\|\le 1$ and $\varphi(\lambda_i)=z_i$ for each
$i$ if and only if the $mn\times mn$ matrix $[(1-z^*_i z_j)
k_{\lambda_i}(\lambda_j)]$ is positive semidefinite. We say  $k$ is a
holomorphic complete interpolation kernel if $k$ is a
holomorphic $m$-interpolation kernel for every $m$.

Our first result is a holomorphic analog of Proposition 2.1.  

\vspace{5mm}

\noindent{\bf Proposition 3.2.}\ \ {\em The kernel $k$ is a holomorphic
$m$-interpolation kernel on $U$ if and only if for each  positive
integer $n$ and each choice of distinct points $\lambda_1, \dots,
\lambda_n\in U$ the kernel $g$ on $N$ defined by $g(i, j) =
k_{\lambda_i} (\lambda_j)$ for $1\le i, j\le n$ is a discrete
$m$-interpolation kernel on $N$.}

\vspace{2mm}

To prove the proposition mimic the argument from [{\bf A2}] that
deduced Theorem 1.27 from Lemma 1.26: Given $z_1, \dots, z_n$ on $\lambda_1,
\dots, \lambda_n$, choose a countable set of uniqueness in
$U$,  and extend $z$ point by point to this set. In the limit, one
gets a bounded operator 
whose adjoint commutes with $M_z^*$ on a dense set, and so comes from a
function $\phi$.

\vspace{2mm}

Now, if $g$ is a discrete kernel on $N$ as in \S2  of this paper,
then there is no particularly distinguished operator of the form
$T_{g, z} $ with $z: N\to \Bbb C$.  However, if $g$ arises as in 
Proposition 3.1 by localizing the holomorphic kernel $k$ to $n$
distinct points $\lambda_1, \dots, \lambda_n$, then Theorem 1.2
provides ample evidence that the $d$-tuple of operators,
$M^*_z | [k_{\lambda_i}, \dots, k_{\lambda_n}]$,  which in the
$T_{g, z}$ notation has the form $T_\lambda = (T_{g, \lambda^1},
\dots, T_{g, \lambda^d})$ where for each $r, \lambda^r: N\to
\Bbb C$ is defined by $\lambda^r(i) = \lambda^r_i$, is highly
distinguished.  To exploit this tuple we introduce the following
notion.

Observe that Theorem 1.2 is equivalent to the following
assertion.  
$$\begin{array}{l}
 {\rm If}\ T\ {\rm has\ an}\ H(\overline K){\rm - extension\ to\ an\
operator\ of\ the\ form}\ \pi(M^*_z)\ {\rm and} \\ 
  T\in S_{\overline \lambda},\ {\rm then}\ T\ {\rm has\ an}\
H(\overline K){\rm - extension\ to\ an\ operator\ of\ the\ form} \\ 
 \pi(M^*_z | [ k_{\lambda_1}, \dots, k_{\lambda_n}]).
\end{array}
\leqno(3.3)
 $$
The assertion of (3.3) is a property of the kernel $k$.  We
introduce the following definition which results from replacing
the 2 occurrences of the word ``extension" by the word
``dilation" in property (3.3).

\vspace{5mm}
\noindent{\bf Definition 3.4.} Let $U\subseteq \Bbb C^d$ be a
bounded domain and let $k$ be a kernel on $U$ of
the type described in the introduction. Then  $k$ is {\it dilation
localizable} if for all choices of $n$ distinct points $\lambda_1,
\dots, \lambda_n \in U$ (3.3) holds with the word extension
replaced by the word  dilation.

We now can state the principal result of this paper.

\vspace{5mm}
\noindent{\bf Theorem 3.5.}  {\em Let $U\subseteq \Bbb C^d$ be  a
bounded domain and let $k$ be a kernel on $U$ of
the type described in the introduction. Then  $k$ is a holomorphic
complete interpolation  kernel if and only if $k$ is  dilation
localizable.}
\vspace{2mm}

We shall deduce Theorem 3.5 as a corollary of the Arveson
dilation machinery  and Theorem 3.6
below.  
If $T\in \cal S_{\lambda}$ it is
clear that  $\sigma(T) = \{\lambda_1, \dots,
\lambda_n\}\subseteq U$ and that the functional calculus map
$$\Phi_T(h) = h(T)\ ,\ \ h\in H(K)$$
 extends by continuity to a continuous unital algebra
homomorphism onto $\cal A_{T}$ of $H^\infty_k$.  Also, it should be
clear on the level of algebra that if $P$ denotes the orthogonal
projection of ${\cal H}$ onto $[k_{\lambda_1}, \dots,
k_{\lambda_n}]$, then
$$\Psi_T(h(PM | [k_{\lambda_1}, \dots, k_{\lambda_n}])) =
h(T)$$
defines a unital algebra homomorphism of $\cal
A_{PM|\, [k_{\lambda_1}, \dots, k_{\lambda_n}]}$ onto
$\cal A_{T}$.  

Recall that $k$ has the
$m$-{\it contractive localization property} if for all positive
integers $n$, all choices of distinct points $\lambda_1, \dots,
\lambda_n \in U$, and all
$T\in \cal S_{\lambda}$ if $\Phi_T$ is $m$-contractive
then $\Psi_T$ is $m$-contractive.  
Define $\rho:H^\infty_k \to \cal A_{PM|\, [k_{\lambda_1}, \dots,
k_{\lambda_n}]}$ by $\rho(h) =h(PM|\, [k_{\lambda_1},
\dots, k_{\lambda_n}])$.  To say that $k$ has the
$m$-contractive localization property means that the diagram
below
can always be completed in the category of operator algebras
with morphisms the unital $m$-contractive algebra
homomorphisms.

\bigskip

\begin{center}
\begin{picture}(100,200)(100,60)
\put(40,150){\vector(1,1){100}}
\put(60,120){\vector(1,-0){180}}
\put(260,150){\vector(1,-1){05}}
\qbezier[30](160,250)(205,205)(257,152)
\put(148,263){\makebox(0,0){ $\cal A_{PM_z |_{ \,
[k_{\lambda_1},\dots, k_{\lambda_n} ]}}$}}
\put(80,210){\makebox(0,0) { $\rho$}}
\put(148,135){\makebox(0,0){ $\Phi_T$}}
\put(30,125){\makebox(0,0) { $H^\infty_k$}}
\put(267,125){\makebox(0,0){ ${\cal A}_{T}$}}
\put(230,194){\makebox(0,0){ $\Psi_T$}}
\end{picture}
\end{center}
\vspace{5mm}

\noindent{\bf Theorem 3.6.}\ \ {\em Let $U\subseteq \Bbb C^d$
be a bounded domain and let $k$ be a kernel on
$U$ of the type described in the introduction. Then $k$ is a
holomorphic
$m$-interpolation kernel if and only if $k$ has the
$m$-contractive localization property.}

\vspace{5mm}

\noindent{\bf Proof.}\ \   Suppose $k$ is an $m$-interpolation
kernel and assume that $T \in S_\lambda$ and 
$\|id_m \otimes \Phi_T \| \le 1$.  We wish
to show that if $h\in \cal M_m(\Bbb C)\otimes H(K)$ and
$\|h(PM_z|\, [k_{\lambda_1}, \dots, k_{\lambda_n} ])\| \le
1$, then $\|h(T)\|\le 1$.  Since $k$ is an $m$-interpolation
kernel, if
$\|h(PM_z|\, [k_{\lambda_1}, \dots, k_{\lambda_n} ])\|\le
1$, then there exists $\tilde h \in \cal M_m(\Bbb C) \otimes
H^\infty_k$ with $\|\tilde h\|_\infty \le 1$ and $\tilde
h(\lambda_i) = h(\lambda_i)$ for each $i$.  

Let us belabor the proof of this last assertion, as it is key.
\begin{eqnarray*}
\Vert h ((PM_z|\, [k_{\lambda_1}, \dots, k_{\lambda_n} ])\|
\ =  \ \| h(PMP) \| \ &=&  \| P h(M) P \| \\
&=& \| P h(M)^* P \| \\
&=& \| \check h (M^*) P \|
\end{eqnarray*}
Now, the fact that $k$ is an $m$-interpolation kernel means that 
$\check h$ can be replaced by ${ h_1}$
where $ h_1 (\overline \lambda_i) = \check h (\overline \lambda_i)$
for each $i$ and 
$$
\| \check h_1 \|_{H^\infty_k} \ = \|  h_1 (M^*) P \| 
= \| \check h (M^*) P \| \leq 1 .
$$
Put $\tilde h = \check h_1$ and it satisfies the assertion.

As $T$ is  in $\cal S_\lambda$, $\tilde h(T)=h(T)$.  Consequently,
$$\|h(T)\| = \|\tilde h(T)\| = \|id_m \otimes \Phi_T (\tilde h)\|
\le 1.$$

Now assume that $k$ has the $m$-contractive localization 
property.  By Proposition 3.2, $k$ will be an $m$-interpolation
kernel if $g_i(j) = k_{\lambda_i} (\lambda_j)$ is a discrete
$m$-interpolation kernel.   We prove this by verifying condition
(iii) of Proposition 2.1.  Accordingly fix $\lambda_1, \dots,
\lambda_{n+1}\in U$ and let $N_1 = \{1, 2, \dots, n+1\}$.  We wish to
show that if $z: N\to \cal M_m(\Bbb C)$, then
$$\|PT_{k, \tilde z}|_{ \cal H_g \ominus \Bbb C g_{n+1}}\| \le \|T_{g, 
z}\| \leqno(3.7)$$
where $P$ denotes the orthogonal projection of $\cal H_g$ onto
${\cal H}_g \ominus \Bbb Cg_{n+1}$ and $\tilde z$ is any extension of  $N$
to $N_1$.  Now exploiting the definition of $g$ it is clear that if the
map $\Omega$ is defined by 
$$\Omega(\check h(M^* |_{ [k_{\lambda_1}, \dots,
k_{\lambda_{n}} ]} )) = \check h (QM^* |_{ [k_{\lambda_1}, \dots,
k_{\lambda_{n+1}} ] \ominus \Bbb Ck_{\lambda_{n+1}}}), \ h \in
H(K)\leqno(3.8)$$
where $Q$  denotes the orthogonal projection onto 
$[k_{\lambda_1}, \dots,
k_{\lambda_{n+1}}] \ominus \Bbb C k_{\lambda_{n+1}}$, then (3.7) is
equivalent to the $m$-contractivity of $\Omega$.  If $\check
\Omega$ is defined by $\check \Omega(S) = \Omega(S^*)^*$ we 
obtain that (3.7) will follow from the $m$-contractivity of
$\check \Omega$.

Now set $T=QM | [k_{\lambda_1}, \dots,
k_{\lambda_{n+1}}]  \ominus \Bbb C k_{\lambda_{n+1}}$.  Evidently,
$\check \Omega = \Psi_T$ and $\Phi_T$ is completely
contractive.  Since in particular $\Phi_T$ is $m$-contractive we
deduce from the $m$-contractive localization property that
$\check \Omega$ is $m$-contractive.  This concludes the proof
of Theorem 3.6.
\ep

{\bf Proof of Theorem 3.5.} By [{\bf Arv1}], the statement that $T^*$ has 
an $H(\overline K)$-dilation to an operator of the form $\pi(M_z^*)$
is equivalent to the assertion that $\Phi_{T^*}$ is a complete contraction;
and that $T^*$ has 
an $H(\overline K)$-dilation to an operator of the form $\pi(M_z^*)
|[k_{\lambda_1},\dots,k_{\lambda_n}]$
is equivalent to the assertion that $\Psi_{T^*}$ is a complete contraction.
So, by interchanging $T$ and $T^*$, the statement of the theorem is
equivalent to saying that $k$ is a holomorphic complete interpolation kernel
if and only if whenever $\Phi_T$ is completely contractive, then $\Psi_T$ is
completely contractive.

It follows from Theorem 3.6 that if $k$ is a holomorphic complete interpolation
kernel, then if $\Phi_T$ is completely contractive, then $\Psi_T$ is
$m$-contractive for all $m$, and hence completely contractive.
Conversely, if the complete contractivity of $\Phi_T$ implies that of
$\Psi_T$,
then the map $\Omega$ of (3.8) will be completely contractive, so the same
argument used in the proof of Theorem 3.6 shows that $k$
is a holomorphic complete interpolation
kernel.

\vspace{1cm}
\noindent {\bf \S 4. Some Examples and Remarks}

\vspace{2mm}
In this section we shall give two concrete applications of
Theorem 3.6.  Our intent  is more to demonstrate that Theorem
3.6 contains interesting function theoretic content  rather than
to work out the most general possible concrete interpolation
theorem that  would follow from the ideas in this section.  Our
first application is a matrix valued generalization of Theorem
0.1 in [{\bf A2}].

\vspace{5mm}

\noindent{\bf Theorem 4.1.}\ \ {\em If $\cal H$ is the Dirichlet
space \Big(i.e. the Hilbert space of analytic functions  on $\Bbb
D$  with reproducing kernel defined by $k_\lambda(\mu) =
\displaystyle\frac{1}{\overline \lambda\mu} \log
\displaystyle\frac{1}{1-\overline \lambda
\mu}\Big)$, then $\cal H$ is a complete interpolation space.}

Our second application constitutes a generalization of the
classical Nevanlinna-Pick interpolation theorem to the ball.

\vspace{5mm}
\noindent{\bf Theorem 4.2.}\ \ {\em Let $d$ be a positive integer
and let $B$ denote the open unit ball in $\Bbb C^d$.  If $\cal H$ is
the Hilbert space of analytic functions on $B$ with reproducing
kernel defined by $k_\lambda(\mu) =\displaystyle
\frac{1}{1-\langle \mu,
\lambda\rangle}$, then $\cal H$ is a complete interpolation
space.}
\vspace{2mm}

Theorem 4.2 should be contrasted with Theorem 3.16 in
 [{\bf A2}]
which derives an analog of the Nevanlinna-Pick interpolation
theorem on the polydisc.  In that solution the first departure of
the canonical interpolation norm from the
$H^\infty$ norm occurs in dimension 3.  In Theorem 4.2 already
when $d=2, H^\infty_k \not= H^\infty (B)$ (this has also been observed in
[{\bf Arv3}]).  Indeed if $n=(n_1,
n_2)$ is a multi-index,  then
$$\|z^n\|^2_{\cal H} = \left(\begin{array}{c}
n_1 +n_2\\ n_2\end{array}\right)^{-1}$$
and
$$\sup\limits_{z\in B}\ |z^n|^2 =
\left(\frac{n_1}{n_1+n_2}\right)^{n_1}
\left(\frac{n_2}{n_1+n_2}\right)^{n_2}.$$
Consequently,
\begin{eqnarray*}
\|z^n\|^2_{H^\infty_k} &\ge& \frac{\|M_{z^n} z^n \|^2_{\cal H} }
{\| z^n \|^2_{\cal H} } \\
&=&
\frac{\|z^{2n}\|^2_{\cal H}}{\| z^n \|^2_{\cal H} }\\
&=& 
\left(\begin{array}{c}
n_1 +n_2\\ n_1\end{array}\right)
\left(\begin{array}{c}
2(n_1 +n_2)\\ 2 n_1\end{array}\right)^{-1}
\left(\frac{n_1}{n_1+n_2}\right)^{-n_1}
\left(\frac{n_2}{n_1+n_2}\right)^{-n_2}
(\sup\limits_{z\in B} |z^n|)^2.
\end{eqnarray*}
Setting $n_1=n_2=j$ and letting $j\to\infty$ we deduce that
the inclusion map of $H^\infty_k$  into $H^\infty (B)$ is not
bounded below.  Hence by the open mapping theorem,
$H^\infty_k \not= H^\infty (B)$.

In connection with the problem of ordinary
$H^\infty$-interpolation on $B$ observe that the  facts that
$H^\infty_k \subseteq H^\infty(B)$ and
$\|\varphi\|_{H^\infty_k} \le \|\varphi\|_{H^\infty(B)}$ whenever
$\varphi \in H^\infty_k$ imply via Theorem 4.2 that
$$\left(\frac{1-z_j\overline z_i}{1-\langle \lambda_j,
\lambda_i\rangle}\right) \ge 0$$
 is a sufficient condition for there to exist a holomorphic
function $\varphi$ on $B$ such that $\varphi(\lambda_i) = z_i$ for
each $i$ and $\sup\limits_{\lambda \in B}\ |\varphi(\lambda)|\le
1$.  On the other hand since $\displaystyle\frac{1}{(1-\langle
\mu,
\lambda\rangle )^d}$ is the Szeg\"o kernel for $B$ it is clear
from the usual proof of the Nevanlinna-Pick on $\Bbb D$ that 
$$\left(\frac{1-z_j\overline z_i}{(1-\langle \lambda_j,
\lambda_i\rangle)^d}\right) \ge 0$$
is a necessary condition for there to exist a holomorphic
function $\varphi$ on $B$ such that $\varphi(\lambda_i)=z_i$ for
each $i$ and $\sup\limits_{\lambda \in B}\ |\varphi(\lambda) | \le
1$.  These considerations prompt one to ask whether there is a
kernel $g$ on the ball (somehow intermediate between
$(1-\langle \mu, \lambda\rangle )^{-1}$ and $(1-\langle \mu,
\lambda\rangle)^{-d})$ with the property that $((1-z_j\overline
z_i) g_{\lambda_i} (\lambda_j)) \ge 0$ is both necessary and
sufficient for ordinary $H^\infty$ interpolation.  We close this
section with an argument which shows that such a kernel $g$
does not exist.

Theorem 4.1 and Theorem 4.2 will both be deduced from the
following fact.

\vspace{5mm}
\noindent{\bf Proposition 4.3.}\ \ {\em Let $\cal H$ be a regular
Hilbert space of analytic functions on a bounded  domain $U\subseteq \Bbb C^d$
with kernel $k$.  If $T\in {\cal F}(k)$
whenever $T$ has a $H(K)$-dilation to an element of ${\cal F}(k)$,
then $\cal H$ is a complete interpolation space.}

\vspace{5mm}
\noindent {\bf Proof.}\ \ The proposition will follow from
Theorem 3.6 if we can establish that $\cal H$ has the
$H(\overline K)$-dilation property.  Accordingly, assume that $\lambda_1,
\dots, \lambda_n$ are $n$  distinct points in $U, T \in \cal
S_{\overline \lambda}$ and $T$ has an $H(\overline K)$-dilation to an
operator of the form $\pi(M^*)$.  We need to  show that $T$ has
an $H(\overline K)$-dilation to an operator of the form $\pi(M^*|\cal
H^\perp_\lambda)$.  This follows immediately from the hypotheses and Theorem 1.2.
\ep

We now prove Theorem 4.1.  Let $\cal H$ denote the Dirichlet
space and set $k_\lambda(\mu)=
\displaystyle\frac{1}{\overline
\lambda\mu}
\log \displaystyle\frac{1}{1-\overline \lambda\mu}$.  By Lemma 2.7 in [{\bf
A2}]
there exists a positive sequence $a_1, a_2, \dots$ such that 
$$\frac{1}{k_\lambda(\mu)} = 1 - \sum\limits^\infty_{n=1}
a_n(\overline \lambda\mu)^n,\leqno(4.4)$$
where the series in (4.4) converges uniformly on compact
subsets of $\Bbb D\times \Bbb D$.  We claim that 
$$T\in \cal F(k)\quad {\rm if\ and \ only \ if}\quad rT\in \cal F(k)\quad
{\rm whenever}\quad 0\le r<1.\leqno(4.5)$$
To prove (4.5) first assume that $T\in \cal F(k)$ and fix $r$ with
$0\le r <1$.  Since $T\in \cal F(k), \sigma(T) \subseteq \Bbb
D^-$.  Hence $\sigma(rT) \subseteq r\Bbb D^-$, a compact
subset of $\Bbb D$.  Thus, by Theorem 1.5 that $rT\in \cal F(k)$
will follow if we can show that $\displaystyle\frac{1}{k} (rT)\ge
0$.  But Theorem 1.3 and  Lemma 1.4 imply that
$\displaystyle\frac{1}{k} (rT)
\ge 0$ follows from the positive definiteness of
$\displaystyle\frac{k_\lambda(\mu)}{k_{r\lambda}(r\mu)}$ on
$\Bbb D$.  Since (4.4) implies that
$$\frac{k_\lambda(\mu)}{k_{r\lambda}(r\mu)} = 1 +
\sum\limits^\infty_{n=1} a_n (1-r^{2n}) (\overline \lambda
\mu)^n k_\lambda(\mu)$$
we deduce that $rT\in \cal F(k)$ as was to be shown.  Now
assume that $rT \in \cal F(k)$ whenever $0\le r <1$.  It follows
from Theorem 1.3 that the hereditary functional calculus map
$$h(M^*) \longrightarrow h(rT)$$
is completely positive whenever $0\le r <1$.  Since for each
$h\in H({\Bbb D}\times \overline {\Bbb D}),\ \lim\limits_{r\to 1^-} h(rT)=h(T)$
we deduce that the map
$$h(M^*) \longrightarrow h(T)$$
is also completely positive.  Hence Theorem 1.3 implies that
$T\in \cal F(k)$.  This establishes (4.5).

The proof of Theorem 4.1 is now easy to conclude.  Assume that
$T\in \cal L(\cal K),$ $ T\in \cal F(k), \cal M\subseteq \cal K$ is a
semi-invariant subspace for $T$, and let  $P$ denote the 
orthogonal projection of $\cal K$ onto $\cal M$.  Theorem 4.1 will
follow from Proposition 4.3 if we can  show that $PT|\cal M\in \cal
F(k)$.  Fix $f\in \cal K$ and let $r<1$.  Employing (4.4) we see that 
\begin{eqnarray*}
\langle \frac{1}{k} (rPT|\cal M)x, x\rangle &=& \|x\|^2 -
\sum\limits^\infty_{n=1} a_nr^{2n} \|PT^nx\|^2\\
\vspace{2mm}
&&\ge \|x\|^2 - \sum\limits^\infty_{n=1} a_nr^{2n} \|T^nx\|^2=
\langle \frac{1}{k} (rT)x, x\rangle.
\end{eqnarray*}
This inequality, (4.5) and Theorem 1.5 imply that
$\displaystyle\frac{1}{k} (rPT|\cal M) \ge 0$.  Since $r<1$ is
arbitrary we conclude via Theorem 1.5 and (4.5) that $PT|\cal
M\in \cal F(k)$ establishing Theorem 4.1.
\ep

The proof of Theorem 4.2 follows the general outline of the
proof of Theorem 4.1.  The equation (4.4) is replaced by the
equation
$$\frac{1}{k_\lambda(\mu)} = 1-\overline \lambda^1 z^1 - \overline
\lambda^2 z^2,$$
 and (4.5) is trivial.

We now prove the claim made at the beginning of this section
that there does not  exist a kernel $k$ on the ball with respect
to which the Nevanlinna-Pick theorem is true in  $H^\infty$
norm.  In fact much more is true.

\vspace{5mm}
\noindent{\bf Proposition 4.6.}\ \ {\em Let $\cal H$ be a regular
Hilbert space of analytic functions on $B$, the open unit ball in
$\Bbb C^d$.  Let $k$ denote the kernel for $\cal H$ and assume
that 0.4 (i) and 0.4 (ii) are  equivalent whenever $m=1, n=2$, and
$\lambda_1, \lambda_2 \in B$.  Then there exists a nonvanishing
holomorphic function $f$ on $B$ such that $k_\lambda(\mu) =
\overline{f(\lambda)} f(\mu) (1-\langle \mu,
\lambda\rangle)^{-1}$.  In particular,  $M$ is unitarily equivalent
to the $d$-tuple $M$ of Theorem 4.2 and if $\lambda_1, \dots,
\lambda_n \in B$ and  $z_1, \dots, z_n \in \cal M_m(\Bbb C)$,
then
$$((1-z_jz^*_i) k_{\lambda_i} (\lambda_j)) \ge 0$$
if and only if
$$\left(\frac{(1-z_jz_i^*)}{1-\langle \lambda_j,
\lambda_i\rangle}\right) \ge 0.$$
}

\vspace{5mm}
\noindent{\bf Proof.}\ \ First observe by considering the
Carath\'eodory metric on $B$ that if $\lambda \in B$ and $z \in
\Bbb D$ then there exists a holomorphic mapping $h:B\to \Bbb D$
with $h(0)=0$ and $h(\lambda)=z$ if  and only if $|z|\le
\|\lambda\|$.  On the other hand the assumed equivalence of 0.4
(i) and 0.4 (ii) implies that there exists a holomorphic mapping
$h:B\to \Bbb D$ with $h(0)=0$ and $h(\lambda)=z$ if and only if 
$$\left[\begin{array}{lc}
k_0(0)& k_0(\lambda)\\
k_\lambda(0)& (1-|z|^2) k_\lambda(\lambda)\end{array}
\right] \ge 0.$$
Hence we conclude that if $\lambda \in B$ and $z\in \Bbb D$, then
$|z| \le \|\lambda\|$ if and only if $|z|^2 k_0(0)
k_\lambda(\lambda) \le k_0(0)
k_\lambda(\lambda)-|k_0(\lambda)|^2$.  In particular,
$$k_\lambda(\lambda) = \frac{|k_0(\lambda)|^2}{k_0(0)}\
\frac{1}{1-|\lambda|^2}\leqno(4.7)$$
whenever $\lambda \in B$.  If we define holomorphic functions $g$
 and $h$ on $B\times B$ by the formulae
$$g(\mu, \lambda) = k_{\overline \lambda}(\mu)$$
and
$$h(\mu, \lambda) = \frac{\overline{k_0(\overline \mu)}
k_0(\lambda)}{k_0(0)}\ \frac{1}{1-\langle \lambda, \overline
\mu\rangle}$$
 we deduce from (4.7) and the fact that  $\{(\mu, \overline \mu):
\mu\in B\}$ is a set of uniqueness that $g = h$, and hence
$k_\lambda(\mu)=\overline{f(\lambda)} f(\mu) (1-\langle \mu,
\lambda\rangle)^{-1} $ if $f$ is defined by
$f(\mu)=k_0(\mu)k_0(0)^{-\frac{1}{2}}$.  The other conclusions of
Proposition 4.6 follow immediately.
\ep

Proposition 4.6 implies that there is no kernel on
the ball with respect to which the Nevanlinna-Pick theorem with
$H^\infty$ norm can be true.  For by Proposition 4.6 such a kernel
could be taken to be defined by $k_\lambda(\mu) = (1-\langle
\mu, \lambda\rangle)^{-1}$.  Since $H^\infty_k \not=
H^\infty(B)$, there exists an $\varphi\in H^\infty_k$ such that
$\|\varphi\|_{H^\infty_k} = \rho >1 = \|\varphi\|_{H^\infty(B)}$. 
To show that the Nevanlinna-Pick theorem doesn't hold, we must find an
$n$-tuple $\lambda_1,\dots,\lambda_n$ such that the matrix
$$
\left( \frac{1 - \varphi(\lambda_j) \overline{\varphi(\lambda_i)}}
{1 - \langle \lambda_j, \lambda_i \rangle} \right)
\leqno(4.8)
$$
is not positive. But the positivity of (4.8) is equivalent to saying
$M_{\varphi}^*$ is a contraction on $[k_{\lambda_1}, \dots, k_{\lambda_n}]$,
and as finite linear combinations of the kernel functions are dense
in $\cal H$, this would contradict the fact that $\|\varphi\|_{H^\infty_k} >
1$.

Note that the kernel in Theorem 4.2 has also been studied by Arveson
[{\bf Arv3}].


\vspace{1cm}

\begin{thebibliography}{Agr-S}

\bibitem[A1]{}
Agler, J., ``The Arveson extension theorem and
coanalytic models", {\it Int. Eq. and Op. Th. }{\bf 5}(1982),
608-631.

 \bibitem[A2] {}Agler, J., ``Some interpolation theorems of
Nevanlinna-Pick type", preprint.

\bibitem[Agr-S]{} Agrawal, O. P. and Salinas, N., ``Sharp kernels
and canonical subspaces", {\it Am. Jour. of Math.} {\bf
110}(1987), 23-48.

\bibitem[Aro]]{} Aronszajn, N., ``Theory of reproducing kernels",
{\it T.A.M.S.} {\bf 68}(1950), 337-404.

\bibitem[Arv1]{} Arveson, W. B., ``Subalgebras of
$C^*$-algebras", {\it Acta Math.} {\bf 123}(1969), 141-224.

\bibitem[Arv2]{}  Arveson, W. B., ``Subalgebras of
$C^*$-algebras II", {\it Acta Math.} {\bf 128}(1972), 271-308.

\bibitem[Arv2]{}  Arveson, W. B., ``Subalgebras of
$C^*$-algebras III: Multivariable Operator Theory'', preprint.

\bibitem[C]{} Curto, R., ``Applications of several complex variables to
multiparameter spectral theory'', {\it Surveys of some recent results in
Operator Theory, Vol. II,} ed. J.B. Conway and B.B. Morrel,
Longman, Harlow, 1988, 25-90.

\bibitem[N]{} Nevanlinna, R., ``\"Uber Beschr\"ankte Funktionen
die in gegebene Punkten vorgeschriebene Werte annehmen",
{\it Ann. Acad. Sci. Fenn. Ser.} {\bf A13}(1919), No. 1.

\bibitem[P]{} Pick, G., ``\"Uber  die  Beschr\"ankungen
analytischer Funktionen, welche durch vorgegebene
Funktionswerte bewirkt werden", {\it Math. Ann.} {\bf
77}(1916),7-23.

\bibitem[T1]{} Taylor, J. L., ``A joint spectrum for several
commuting operators", {\it Jour. Func. Anal.} {\bf 6}(1970),
172-191.

\bibitem[T2]{} Taylor, J. L., ``The analytic functional calculus
for several commuting operators", {\it Acta Math.} {\bf
125}(1970), 1-38.


\end{thebibliography}
\end{document}